\title{\vspace{-0.1cm} Pancyclicity of Hamiltonian and highly connected graphs}
\author{
Peter Keevash \thanks{School of Mathematical Sciences,
Queen Mary, University of London, Mile End Road, London E1 4NS, UK.
Email: p.keevash@qmul.ac.uk.
Research supported in part by NSF grant DMS-0555755.}
\and
Benny Sudakov\thanks{Department of Mathematics,
UCLA, Los Angeles, 90095. E-mail: bsudakov@math.ucla.edu. Research
supported in part by NSF CAREER award DMS-0812005
and a USA-Israeli BSF grant.}
}
\documentclass[11pt]{article}
\usepackage{amsfonts,amssymb,amsmath,latexsym,graphicx}

\def\qed{\ifvmode\mbox{ }\else\unskip\fi\hskip 1em plus 10fill$\Box$}

\newtheorem{theo}{Theorem}[section]
\newtheorem{lemma}[theo]{Lemma}
\newtheorem{prop}[theo]{Proposition}

\newcommand{\nib}[1]{\noindent {\bf #1}}

\newcommand{\sm}{\setminus}

\newcommand{\sub}{\subseteq}

\begin{document}
\date{}

\maketitle

\begin{abstract}
A graph $G$ on $n$ vertices is {\em Hamiltonian} if it contains a cycle of length $n$
and {\em pancyclic} if it contains cycles of length $\ell$ for all $3 \le \ell \le n$.
Write $\alpha(G)$ for the {\em independence number} of $G$, i.e.\ the size of the largest
subset of the vertex set that does not contain an edge, and $\kappa(G)$ for the (vertex)
{\em connectivity}, i.e.\ the size of the smallest subset of the vertex set that can be deleted
to obtain a disconnected graph. A celebrated theorem of Chv\'atal and Erd\H{o}s says that
$G$ is Hamiltonian if $\kappa(G) \ge \alpha(G)$. Moreover, Bondy suggested that almost any non-trivial
conditions for Hamiltonicity of a graph should also imply pancyclicity. Motivated by this,
we prove that if $\kappa(G) \ge 600\alpha(G)$
then $G$ is pancyclic. This establishes a conjecture of Jackson and Ordaz up to a constant factor.
Moreover, we obtain the more general result that if $G$ is Hamiltonian with
minimum degree $\delta(G) \ge 600\alpha(G)$ then $G$ is pancyclic. Improving an old result of
Erd\H{o}s, we also show that $G$ is pancyclic if it is Hamiltonian  and $n \ge 150\alpha(G)^3$.
Our arguments use the following theorem of independent interest
on cycle lengths in graphs: if $\delta(G) \ge 300\alpha(G)$ then
$G$ contains a cycle of length $\ell$ for all $3 \le \ell \le \delta(G)/81$.
\end{abstract}

\section{Introduction}

A {\em Hamilton cycle} is a spanning cycle in a graph, i.e.\ a cycle passing through all vertices.
A graph is called {\em Hamiltonian} if it contains such a cycle. Hamiltonicity is one of the most
fundamental notions in graph theory, tracing its origins to Sir William Rowan Hamilton in the 1850's.
Deciding whether a given graph contains a Hamilton cycle is NP-complete,
so we do not expect to have a simple characterisation for this property.
There is a vast literature in graph theory devoted to obtaining sufficient conditions
for Hamiltonicity. For more details we refer the interested reader
to the surveys of Gould \cite{G1,G2}. The classical result giving such a condition is
Dirac's theorem \cite{D}, which says that every graph $G$ with $n \geq 3$ vertices
and minimum degree at least $n/2$ is Hamiltonian.
This theorem was generalised by Bondy \cite{B1}, who showed that the same assumptions imply
that $G$ is {\em pancyclic}, i.e.\ contains cycles of length $\ell$ for all $3 \le \ell \le n$.
In \cite{B2} Bondy proposed the following `metaconjecture' which has had a considerable
influence on research on cycles in graphs.

\nib{Metaconjecture.} Almost any non-trivial condition on a graph which implies that the
graph is Hamiltonian also implies that the graph is pancyclic.
(There may be a simple family of exceptional graphs.)

Another classical condition for a graph to be Hamiltonian
is given by a theorem of Chv\'atal and Erd\H{o}s \cite{CE},
who showed that if a graph $G$ satisfies $\kappa(G) \ge \alpha(G)$ then it is Hamiltonian.
Here $\alpha(G)$ is the {\em independence number}, i.e.\ the size of the largest subset of the
vertex set that does not contain an edge, and $\kappa(G)$ is the (vertex) {\em connectivity},
i.e.\ the size of the smallest subset of the vertex set that can be deleted to obtain a disconnected
graph. Motivated by Bondy's metaconjecture, Amar, Fournier and Germa \cite{AFG}
obtained several results on the lengths of cycles in a graph $G$ that satisfies the
Chv\'atal-Erd\H{o}s condition $\kappa(G) \ge \alpha(G)$. They conjectured that if such a graph $G$
is not bipartite then $G$ contains cycles of length $\ell$ for all $4 \le \ell \le n$.
(The case when $G=C_5$ is a $5$-cycle needs to be excluded.)
Note that the balanced complete bipartite graph $K_{k,k}$ satisfies $\kappa(G)=\alpha(G)=k$
but is not pancyclic, indeed it has no odd cycles. They also made the weaker conjecture that
the same conclusion holds under the additional assumption that $G$ is triangle-free.
Lou \cite{L} proved the stronger statement that if $\kappa(G) \ge \alpha(G)$ and $G$ is
triangle-free then $G$ contains cycles of length $\ell$ for all $4 \le \ell \le n$, unless
$G=C_5$ or $G=K_{k,k}$ for some $k$. Note that the connectivity $\kappa(G)$ is bounded above by the
degree of any vertex. If $\kappa(G) > \alpha(G)$ then there is an edge inside any neighbourhood of $G$,
so in particular $G$ must contain triangles. Jackson and Ordaz \cite{JO} conjectured that
if $\kappa(G)> \alpha(G)$ then $G$ is pancyclic.

To approach these conjectures it is natural to try to prove pancyclicity under a
stronger connectivity assumption. A remarkable theorem of Erd\H{o}s \cite{E},
proving a conjecture of Zarins,
shows that instead of making a connectivity assumption,
it suffices to assume that $G$ is Hamiltonian and the number of vertices
is sufficiently large compared to the independence number.
He showed that if $G$ is a Hamiltonian graph on $n$ vertices with $\alpha(G)=k$
and $n>4k^4$ then $G$ is pancyclic. It then follows from \cite{CE} that
$\kappa(G)\ge 4(\alpha(G)+1)^4$ is sufficient for pancyclicity.
(Various considerably weaker results were subsequently obtained by Flandrin et al.,
e.g.\ \cite{FLMW}, who were presumably unaware of Erd\H{o}s' paper.)
Our main result improves this bound significantly: we show that a connectivity which is only linear in
the independence number suffices for pancyclicity. This establishes the conjecture of
Jackson and Ordaz mentioned above up to a constant factor.
Moreover, we prove that pancyclicity already follows from assuming that $G$ is Hamiltonian
with minimum degree $\delta(G)$ at least linear in the independence number.

\begin{theo}\label{pan-mindeg}
If $G$ is a Hamiltonian graph with $\delta(G) \ge 600\alpha(G)$ then $G$ is pancyclic.
In particular, if $G$ is any graph with $\kappa(G) \ge 600\alpha(G)$ then $G$ is pancyclic.
\end{theo}

Erd\H{o}s \cite{E} remarked that the bound $n>4k^4$ in his result is unlikely to be best possible. He also
noticed that a quadratic lower bound for $n$ in terms of $k$ is necessary
for Hamiltonicity to imply pancyclicity. Our next theorem improves Erd\H{o}s' result and
shows that a cubic dependence of $n$ on $k$ is already sufficient.

\begin{theo}\label{pan-n}
If $G$ is a Hamiltonian graph with $|V(G)| \ge 150\alpha(G)^3$ then $G$ is pancyclic.
\end{theo}

Our arguments will use a theorem of independent interest on cycle lengths in graphs.

\begin{theo} \label{short-cycles}
If $G$ is a graph with $\delta(G) \ge 300\alpha(G)$
then $G$ contains a cycle of length $\ell$ for all $3 \le \ell \le \delta(G)/81$.
\end{theo}

It is instructive to compare Theorem \ref{short-cycles} with a result of
Nikiforov and Schelp \cite{NS}, who showed that when the minimum degree $\delta(G)$ is linear
in the number of vertices then $G$ contains even cycles of all lengths between $4$ and $\delta(G)+1$
and, after excluding some exceptional cases, odd cycles of all lengths between a constant and $\delta(G)+1$.
We refer the reader to the chapter of Bondy in \cite{B3} for other results on cycle lengths
in graphs, and to \cite{GHS,V,SV} as examples of more recent related results.

Next we describe a simple example showing that Theorems \ref{pan-mindeg} and \ref{short-cycles}
are best possible up to the constant factors and that the lower bound for $|V(G)|$ in
Theorem \ref{pan-n} has to be at least quadratic in $\alpha(G)$. Suppose $k \ge 3$ and
let $G$ be the graph on $n=k(2k-2)$ vertices obtained by taking $k$ vertex-disjoint cliques
$X_1,\cdots,X_k$ of size $2k-2$ and adding a matching of size $k$ which has exactly one edge between
$X_i$ and $X_{i+1}$ for all $1 \le i \le k$ (where $X_{k+1}:=X_1$).
Then it is easy to check that $G$ is Hamiltonian, $\alpha(G)=k$ and $\delta(G)=2k-3$,
but $G$ does not contain a cycle of length $2k-1$, so is not pancyclic.

The organisation of this paper is as follows. In the next section we collect various
known results that we will use in our arguments. We present the proofs of our theorems
in Section 3. The final section contains some concluding remarks.
We systematically omit rounding signs for the sake of clarity of presentation.
We also do not make any serious attempt to optimise
absolute constants in our statements and proofs.

\nib{Notation.}
Suppose $G$ is a graph. For a vertex $v$ we let $N(v)$ denote its neighbourhood and
$d(v)=|N(v)|$ its degree. If $X$ is a set of vertices then $G[X]$ is the restriction
of $G$ to $X$, i.e.\ the graph with vertex set $X$ whose edges are edges of $G$ with
both endpoints in $X$. We write $e_G(X)=e(G[X])$ for the number of edges in
$X$. If $X$ and $Y$ are sets of vertices then $e_G(X,Y)$ is the number of edges with
one endpoint in $X$ and the other in $Y$. We omit the subscript $G$ if there is no
danger of ambiguity. A walk in $G$ is a sequence of vertices $W = x_0 \cdots x_t$
such that $x_i$ is adjacent to $x_{i+1}$ for $0 \le i \le t-1$. (The vertices need
not be distinct.) The length $\ell(W)=t$ of $W$ is the number of edges in $W$,
counting multiplicity of repeated edges. A path is a walk in which no vertices are repeated.
A cycle is a walk in which no vertices are repeated, except that the first and last vertices are equal.

\section{Preliminaries}

In this section we collect various results that will be used in our arguments.
We include the short proofs for the convenience of the reader.

\subsection{Degrees}

We start with two well-known propositions on degrees in graphs.

\begin{prop}\label{bip}
Suppose $G$ is a graph with minimum degree at least $d$.
Then $G$ has a bipartite subgraph $B$ with minimum degree at least $d/2$.
\end{prop}

\nib{Proof.} Consider a bipartite subgraph $B$ of $G$ with as many edges as possible.
Let $X$ and $Y$ be the two parts of the bipartition of $B$.
Then any vertex $v \in X$ has at least $d(v)/2$ neighbours in $Y$,
or we could improve the partition by moving $v$ to $Y$.
The same argument applies for $v \in Y$. \qed

\begin{prop}\label{avmindeg}
Suppose $G$ is a graph with average degree at least $d$.
Then $G$ has an induced subgraph with minimum degree at least $d/2$.
\end{prop}

\nib{Proof.} Suppose $G$ has $n$ vertices and construct a sequence of graphs
$G_n=G, G_{n-1}, \cdots$ where if $G_i$ has minimum degree less than $d/2$ we
construct $G_{i+1}$ by deleting a vertex of degree less than $d/2$ from $G_i$.
The number of edges deleted in this process is less than $nd/2 \le e(G)$
so it must terminate at some induced subgraph with minimum degree at least $d/2$. \qed

\subsection{Breadth first search trees}

Suppose that $B$ is a graph and $x$ is a vertex of $B$.
We construct a {\em breadth first search tree} $T$ in $B$ starting
at $x$ by the following iterative procedure.
We start with $T_0$ equal to the one-vertex tree on $x$.
Then at step $i \ge 1$, we let $N_i$ be the set of vertices not in $T_{i-1}$ that
have at least one neighbour in the tree $T_{i-1}$, and construct
$T_i$ on the vertex set $V(T_{i-1}) \cup N_i$ by adding an edge
from each vertex $v$ in $N_i$ to a neighbour of $v$ in $T_{i-1}$.

\begin{prop}\label{bfs}
Suppose $B$ is a bipartite graph and $T$ is a breadth first search
tree in $B$ starting from a vertex $x$. Let $N_i$ be the set of
vertices at distance $i$ from $x$ in $T$. Then any vertex in $N_i$
is at distance $i$ from $x$ in $B$. Also, each $N_i$ is an independent
set in $B$ and all edges of $B$ join $N_i$ to $N_{i+1}$ for some $i \ge 0$.

Now suppose also that $B$ has $n$ vertices and minimum degree $d \ge 5$.
Then there is some number $i \ge 0$ such that
$e_B(N_i,N_{i+1}) \ge \frac{d}{4}(|N_i|+|N_{i+1}|)$.
Furthermore, if $m \ge 0$ is the smallest number with
$e_B(N_m,N_{m+1}) \ge \frac{2d}{9}(|N_m|+|N_{m+1}|)$
then $m \ge 1$ and $|N_{i+1}| \ge 2|N_i|$ for $0 \le i \le m-1$.
\end{prop}

\nib{Proof.} Let $T_0, T_1, \cdots$ be the sequence of trees in the
breadth first search construction. For any $v$ in $N_i$, the neighbours in $B$
of $v$ within $V(T_{i-1})$ must lie in $N_{i-1}$, or we would have already added $v$
to $T_{i-1}$. We deduce that the distance in $B$ from $v$ to $x$ is $i$.
Next consider any $y$ and $z$ in $N_i$, let $P$ be the path between $y$ and $z$
in $T_i$ and let $w$ be the closest point to $x$ on $P$. If $w \in N_j$ then
the length of $P$ is $2(i-j)$, which is even, so $yz$ cannot be an edge,
since $B$ is bipartite. This shows that $N_i$ is independent and
all edges of $B$ join $N_i$ to $N_{i+1}$ for some $i \ge 0$.

Now suppose that $B$ has $n$ vertices and minimum degree at least $d$, so that $e(B) \ge dn/2$.
We cannot have $e_B(N_i,N_{i+1}) < \frac{d}{4}(|N_i|+|N_{i+1}|)$ for all $i \ge 0$,
as this would give the contradiction
$$dn/2 \le e(B) = \sum_{i \ge 0} e(N_i,N_{i+1}) < \frac{d}{4} \sum_{i \ge 0} (|N_i|+|N_{i+1}|)
= \frac{d}{4}(2n-1).$$
Consider the smallest $m \ge 0$ with $e_B(N_m,N_{m+1}) \ge \frac{2d}{9}(|N_m|+|N_{m+1}|)$.
Then $m \ge 1$, since $e_B(N_0,N_1)=|N_1|$ and $d \ge 5$,
so $\frac{2d}{9}(|N_0|+|N_{1}|) \ge |N_1|+1$.
We claim that $|N_{i+1}| \ge 2|N_i|$ for $0 \le i \le m-1$.
For suppose that this is not the case, and consider the smallest
$0 \leq i \leq m-1$ for which $|N_{i+1}| < 2|N_i|$. Then $i \ge 1$, since
$|N_0|=1$ and $|N_1| \ge d \ge 5$.
There are at least $d|N_i|$ edges incident to $N_i$,
so we must have $e_B(N_{i-1},N_i) \ge d|N_i|/3$ or $e_B(N_i,N_{i+1}) \ge 2d|N_i|/3$.
In the first case we have
$e_B(N_{i-1},N_i) \ge d|N_i|/3 \ge \frac{2d}{9}(|N_{i-1}|+|N_i|)$,
since $|N_i| \ge 2|N_{i-1}|$. In the other case we have
$e_B(N_i,N_{i+1}) \ge 2d|N_i|/3 \ge \frac{2d}{9}(|N_i|+|N_{i+1}|)$,
since $|N_{i+1}| < 2|N_i|$. Either way we have a contradiction to
the minimality of $m$, so the claim is proved. \qed

\subsection{Independence number}

Here we give some well-known relationships between degrees, chromatic number
and independence number.

\begin{prop} \label{indep}
Suppose $G$ is a graph on $n$ vertices with maximum degree at most $k$.
Then $G$ contains an independent set of size at least $n/(k+1)$.
\end{prop}

\nib{Proof.} Construct an independent set $S$ greedily by repeatedly choosing
any currently available vertex and then marking all its neighbours as unavailable.
At the end of this process every vertex of $G$ is either in $S$ or marked as
unavailable. At most $k|S|$ vertices have been marked unavailable,
so $n \le |S|+k|S|$, i.e.\ $|S| \ge n/(k+1)$. \qed

\begin{prop} \label{chromatic}
Suppose $G$ is a graph for which every induced subgraph has a vertex of degree at most $k$.
Then $G$ has chromatic number at most $k+1$.
\end{prop}

\nib{Proof.}
Define a sequence of induced subgraphs $G_n,\cdots,G_0$ starting from $G_n=G$,
where $G_{i-1}$ is obtained from $G_i$ by deleting a vertex $v_i$ of degree at most $k$.
Consider the vertices in the order $v_1,\cdots,v_n$ and greedily colour them using
$\{1,\cdots,k+1\}$. When we colour $v_i$ we have used at most $k$ colours on its neighbours
in $G_i$, so there is an available colour in $\{1,\cdots,k+1\}$. \qed

\begin{prop} \label{indep-mindeg}
Suppose $G$ is a graph with independence number $\alpha(G) \le k$
and $n \ge dk+1$ vertices. Then $G$ contains an induced subgraph $H$ with at most $dk+1$
vertices and minimum degree at least $d$.
\end{prop}

\nib{Proof.} Let $S$ be a set of $dk+1$ vertices of $G$.
The restriction $G[S]$ of $G$ to $S$ must have chromatic number at least $d+1$,
otherwise it would contain an independent set of size at least $|S|/d > k$,
contradicting our assumption on $G$. By Proposition \ref{chromatic} $G[S]$
contains an induced subgraph $H$ with minimum degree at least $d$. \qed

\subsection{Paths}

Now we give some simple tools for manipulating the lengths of paths
in cycles when there is a bound on the independence number.

Given any path or cycle $W$ in a graph $G$, we refer to any set of $t+1$ consecutive points
on $W$ as a {\em $t$-interval} (so $t$ is the length of the interval).
If $J$ is an interval of length at least $2$ such that the endpoints of $J$ are adjacent
then we call $J$ a {\em jump} of $W$ in $G$.
For a jump $J$ we write $\partial J$ for the edge joining the ends of $J$
and $J^o$ for the subinterval of internal points obtained by removing its ends.

\begin{prop}\label{jump}
Suppose $G$ is a graph with independence number $\alpha(G) \le k$,
$W$ is a path or cycle in $G$ and $I$ is an interval of length at least $2k$ on $W$.
Then $I$ contains a jump of $W$ in $G$ with length at most $2k$.
\end{prop}

\nib{Proof.} Starting at one end of $I$ consider the points with positions
$1,3,5,\cdots,2k+1$. This set of $k+1$ points must contain an edge, since
$\alpha(G) \le k$. \qed

\begin{prop}\label{shorten}
Suppose $G$ is a graph with independence number $\alpha(G) \le k$
and $P$ is a path of length $p$ in $G$ joining two vertices $x$ and $y$.
Then for any number $1 \leq q \leq p$ there is a path of some length $\ell$
joining $x$ and $y$ with $q \leq \ell \leq q+2k-2$.
\end{prop}

\nib{Proof.} We use induction on $p$. The statement is clearly true for $p \leq 2k-1$,
so suppose that $p \geq 2k$. By Proposition \ref{jump} there is a jump $J$ of $P$ in $G$
of some length $j$ with $2 \le j \le 2k$. Replacing the portion of $P$ along $J$ by the
edge $\partial J$ joining the ends of $J$ gives a path $P'$ of length $p-j+1$ joining $x$ and $y$.
Now for all $q>p-j+1$ we can use the original path $P$,
and for all $q \leq p-j+1$ we can apply the induction hypothesis to $P'$. \qed

We also need the following well-known proposition.

\begin{prop}\label{length}
Suppose $G$ is a graph with minimum degree at least $d$
and $x$ is a vertex of $G$.
Then $G$ contains a path of length at least $d$ starting at $x$.
Furthermore, if $G$ is bipartite then $G$ contains such
a path of length at least $2d-1$.
\end{prop}

\nib{Proof.} Let $P$ be a longest path in $G$ starting at $x$
and let $y$ be the last vertex of $P$. By the minimum degree condition $y$
has at least $d$ neighbours, and these all belong to $P$ by choice of a
longest path, so $P$ contains at least $d+1$ vertices.
Furthermore, if $G$ is bipartite, then $y$ is not adjacent to
any vertex at even distance from $y$ along $P$,
so $P$ contains at least $2d$ vertices. \qed

\subsection{Hamiltonicity}

Here we give two more substantial lemmas on Hamiltonian graphs which appeared implicitly in \cite{E}.
One facilitates absorption of a vertex to create a Hamiltonian graph
with one more vertex, the other deletion of a vertex to create a Hamiltonian
graph with one fewer vertex.

\begin{lemma} \label{add}
Suppose $G$ is a graph, $x$ is a vertex of degree at least $k+1$ in $G$ and
$H = G \sm \{x\}$ is a Hamiltonian graph with independence number $\alpha(H) \le k$.
Then $G$ is Hamiltonian.
\end{lemma}

\nib{Proof.}
Suppose that $H$ has $n$ vertices. Label them with $[n]=\{1,\cdots,n\}$
such that $\{i,i+1\}$ is an edge for $1 \le i \le n$, where addition
is mod $n$, i.e.\ $n+1$ is identified with $1$. Let $A \sub [n]$ be the
neighbourhood of $x$ and let $A^+ = \{a+1: a \in A\}$.
Since $|A^+| \ge d(x) \ge k+1 > \alpha(H)$ there is an edge
$\{y,z\}$ in $A^+$, where without loss of generality $y<z$.
Now we can form a Hamilton cycle in $G$ by starting at $x$,
going to $z-1 \in A$, decreasing to $y$, using the edge $\{y,z\}$
to get to $z$, increasing to $n$, going to $1$, increasing to $y-1 \in A$,
then ending at $x$. \qed

We remark that the argument in Lemma \ref{add} is the main idea
in the proof of the Chv\'atal-Erd\H{o}s theorem.

\begin{lemma}\label{erdos}
Suppose $G$ is a Hamiltonian graph on $n \ge (2k+1)(k^2+k+1)$ vertices with
independence number $\alpha(G) \le k$. Then $G$ contains a cycle of length $n-1$.
\end{lemma}

\nib{Proof.} Choose a Hamilton cycle $C$ in $G$ and label the vertices as
$v_1,\cdots,v_n$ so that the edges of $C$ are $v_iv_{i+1}$ for $1 \le i \le n$.
(As above we use the convention $v_{n+1}=v_1$.)
Set $s=k^2+k+1$ and let $I_1,\cdots,I_s$ be disjoint $2k$-intervals in $C$.
Proposition \ref{jump} gives jumps $J_1,\cdots,J_s$,
where each $J_i$ is a subinterval of $I_i$ of length at least $2$
and the ends of $J_i$ are adjacent in $G$. We say that $J_i$ is {\em good}
if each internal vertex $v \in J_i^o$ has at least $k+1$ neighbours in $V(G) \sm J_i^o$.
We claim that there is a good jump. For suppose to the contrary that we can
choose $v_i \in J_i^o$ such that $v_i$ has at most $k$ neighbours in $V(G) \sm J_i^o$
for $1 \le i \le s$. Then $\{v_1,\cdots,v_s\}$ spans a subgraph of $G$ with
maximum degree at most $k$, so by Proposition \ref{indep} contains an independent
set of size bigger than $k$, contradicting our assumption on $G$.
Thus there is a good jump, say $J_1$. Now we construct a cycle of length $n-1$ as follows.
First we replace the portion of $C$ traversing the jump $J_1$ with the edge $\partial J_1$
between the endpoints of $J_1$. Then we use Lemma \ref{add} to put back the vertices
of $J_1^o$ one by one, increasing the length of the cycle until only one vertex has
not been replaced. \qed

\subsection{Cycles}

The following lemma of Erd\H{o}s, Faudree, Rousseau and Schelp \cite{EFRS}
will allow us to find a cycle of some particular length,
using a breadth first search tree and the independence assumption.
The proof of this lemma can be found in rather abbreviated form within
the proof of Theorem 1 in \cite{EFRS}. For the convenience of the reader
we include a proof here.

\begin{lemma} \label{efrs}
Suppose $G$ is a graph containing no cycle of length $\ell$,
$T$ is a tree in $G$, $v$ is a vertex of $T$, $h < \ell/2$,
and $Z$ is the set of vertices at distance $h$ in $T$ from $v$.
Then the restriction of $G$ to $Z$ is $(\ell-2)$-colourable,
and so $Z$ contains a subset of size at least
$|Z|/(\ell-2)$ that is independent in $G$.
\end{lemma}

\nib{Proof.}
Fix a plane drawing of $T$ such that in the $(x,y)$-coordinate system
$v$ is at the origin and for $i \ge 0$ points at distance $i$ from $v$ have $x$-coordinate $i$.
Say that a path $z_0z_1 \cdots z_t$ in $G$ using vertices of $Z$
is increasing if the $y$-coordinates of the vertices $z_0,z_1,\cdots,z_t$
forms an increasing sequence.

Given any increasing path $P=z_0z_1 \cdots z_t$ we let $P'$ be the unique path in $T$
from $z_0$ to $z_t$ and we let $v_P$ be the closest point to $v$ on $P'$.
We claim that we can remove either $z_0$ or $z_t$ to obtain a path $Q$
such that $v_Q=v_P$ and $Q'$ has the same length as $P'$.
To see this, we observe that it can only fail if
the path in $T$ from $z_1$ to $v_P$ meets the path
from $z_t$ to $v_P$ before it reaches $v_P$
and the path in $T$ from $z_{t-1}$ to $v_P$ meets the path
from $z_0$ to $v_P$ before it reaches $v_P$.
But this would contradict our choice of a plane drawing of $T$, so the claim holds.

Now we claim that there is no increasing path of length $\ell-2$.
For suppose that $P=z_0z_1 \cdots z_{\ell-2}$ is an increasing path.
Let $\ell'$ be the length of the path $P'$ in $T$ from $z_0$ to $z_{\ell-2}$.
Then $2 \le \ell' \le 2h \le \ell-1$.
We construct a sequence of paths $P_0=P, P_1, \cdots, P_{\ell-3}$ where each
$P_{i+1}$ is obtained from $P_i$ by removing an endpoint in such a way
that $v_{P_{i+1}}=v_{P_i}$, so that for each $i$ we have $v_{P_i}=v_P$
and $P'_i$ has length $\ell'$. Since $P_i$ has length $\ell-2-i$,
$P_i \cup P'_i$ forms a cycle of length $\ell'+\ell-2-i$.
Setting $i=\ell'-2$ we obtain a cycle of length $\ell$,
which contradicts our assumption on $G$, so the claim holds.

Finally, we define a colouring $c:Z\to\{0,1,\cdots,\ell-3\}$
where $c(z)$ is the length of the longest increasing path starting at $z$.
This is a proper colouring of $G[Z]$, as if $z,z' \in Z$ with $z$ below $z'$ (say)
then we can add $zz'$ to any increasing path starting at $z'$,
so $c(z)>c(z')$. Since all colour classes of $c$ are independent
we have an independent set of size at least $|Z|/(\ell-2)$. \qed

\subsection{Probability}

Finally we record the standard Chernoff bounds for large deviations
of binomial random variables.

\begin{lemma} (Chernoff bounds, see \cite{AS} Appendix A)
Suppose $X$ is a binomial random variable with parameters $(n,p)$ and $a \ge 0$.
\begin{itemize}
\item[(i)] If $p=1/2$ then $\mathbb{P}(X-n/2>a)
= \mathbb{P}(X-n/2<-a) < e^{-2a^2/n}$.
\item[(ii)] $\mathbb{P}(X-np>a) < e^{-a^2/2pn+a^3/2(pn)^2}$.
\item[(iii)] $\mathbb{P}(X-np<-a) < e^{-a^2/2pn}$.
\end{itemize}
\end{lemma}

\section{Proofs}

In this section we present proofs of our three theorems. Throughout we will suppose that $G$ is
a graph with independence number $\alpha(G) \le k$. Also, we can suppose $k \ge 2$,
otherwise we have the trivial case when $G$ is a complete graph.
We start with a lemma that provides two vertices that are connected
by paths with every length in some interval.

\begin{lemma} \label{consecutive-paths}
Suppose $G$ is a graph on $n$ vertices with independence number $\alpha(G) \le k$
and $B$ is a bipartite subgraph of $G$ with minimum degree $\delta(B) = d > 9k/2$.
Suppose $x$ is a vertex of $B$ and let $N_i$ denote the set of vertices
at distance $i$ from $x$ in $B$. Let $m \ge 1$ be the smallest number with
$e_B(N_m,N_{m+1}) \ge \frac{2d}{9}(|N_m|+|N_{m+1}|)$. Then
\begin{itemize}
\item[(i)] $|N_m| \ge 2^{m-1}d$, $m \le \log_2 \left( \frac{n+d-1}{d} \right)$
and $G$ contains cycles of length $\ell$ for all $3 \le \ell \le |N_m|/k$,
\item[(ii)] there are sets $N'_m \sub N_m$ and $N'_{m+1} \sub N_{m+1}$
forming the parts of a bipartite subgraph $B'$ of $B$ with minimum degree at least $2d/9$,
\item[(iii)] there is a vertex $y$ in $B'$ such that there is a path between
$x$ and $y$ in $G$ of length $\ell$, for any $\ell$ with $m \le \ell \le m+4d/9-2$.
\end{itemize}
\end{lemma}

\nib{Proof.}
By Proposition \ref{bfs} we have $|N_{i+1}| \ge 2|N_i|$ for $0 \le i \le m-1$.
Since $|N_1| \ge \delta(B) = d$ we have $|N_i| \ge 2^{i-1}d$ for $1 \le i \ge m$.
Applying Lemma \ref{efrs} to $Z=N_i$ for $1 \le i \le m$
we see that $G$ contains cycles of length $\ell$ for $2i+1 \le \ell \le |N_i|/k$.
Since $d > 9k/2$ and $|N_i| \ge 2^{i-1}d$, it is easy to check that $2i+2 \leq |N_i|/k$,
so the intervals $[2i+1, |N_i|/k], 1 \leq i \leq m$ together contain
all integers from $3$ to $|N_m|/k$. Also, $n \ge \sum_{i=0}^m |N_i| \ge 1+(2^m-1)d$
gives the required bound on $m$, so statement (i) holds.
Statement (ii) follows from Proposition \ref{avmindeg}. Indeed,
$N_m$ and $N_{m+1}$ form the parts of a bipartite
subgraph of $B$ with average degree at least $4d/9$. Thus it contains a subgraph $B'$ with parts
$N'_m \sub N_m$ and $N'_{m+1} \sub N_{m+1}$ with minimum degree at least $d'=2d/9$.
Since $d > 9k/2$ the minimum degree in $B'$ is at least $k+1$. In particular,
$|N'_m| \ge k+1$. Since $\alpha(G) \le k$ there is an edge $yz$ of $G$ in $N'_m$.
We claim that this choice of $y$ satisfies statement (iii).
To see this we give separate arguments for paths of length $m+2t$, $t \ge 0$
and paths of length $m+2t+1$, $t \ge 0$. By Proposition \ref{length},
for $0 \le 2t \le 2d'-2=4d/9-2$ there is a path of length $2t$ in $B'$ from $y$
to a vertex $w$ in $N'_m$, which can be combined with the path in $T$
from $w$ to $x$ to give a path of length $m+2t$ between $x$ and $y$.
Next, consider the bipartite graph $B' \sm \{y\}$, which has minimum degree at least $d'-1$.
Then, again by Proposition \ref{length}, for $0 \le 2t \le 2d'-4$
we can find a path in $B'$ of length $2t$ from $z$ to a vertex $w \in N'_m$,
which can be combined with the edge $yz$ and the path in $T$
from $w$ to $x$ to give a path of length $m+2t+1$ between $x$ and $y$. \qed

Now we prove our first result, which states that
a graph $G$ on $n$ vertices with independence number
$\alpha(G) \le k$ and minimum degree $\delta(G) = d \ge 300k$
contains a cycle of length $\ell$ for all $3 \le \ell \le d/81$.

\noindent
{\bf Proof of Theorem \ref{short-cycles}.}\,
By Proposition \ref{bip} we can choose a bipartite subgraph $B$ of $G$
with minimum degree $\delta(B) \ge d/2$. Fix any vertex $x$ and let $N_i$
be the set of vertices at distance $i$ from $x$ in $B$.
By Lemma \ref{consecutive-paths}, for some $m \geq 1$ we have
cycles in $G$ of length $\ell$ for all $3 \le \ell \le |N_m|/k$,
where $|N_m| \ge 2^{m-1}(d/2)=2^{m-2}d$. We also have subsets $N'_m \sub N_m$ and $N'_{m+1} \sub N_{m+1}$
spanning a bipartite subgraph $B'$ of $B$ with minimum degree at least $\frac{2}{9}(d/2)=d/9$.
We can assume that $|N_m| < kd/81$, since otherwise we are done.
Also, by choosing $d/9$ neighbours in $N'_{m+1}$ for each vertex in $N'_m$
and deleting all other vertices of $N'_{m+1}$ we can assume that $|N'_{m+1}| < kd^2/729$.
Next we consider a partition $N'_m = P \cup Q$, where each vertex of $N'_m$
is randomly and independently placed in $P$ or $Q$ with probability $1/2$.
Since every vertex in $N'_{m+1}$ has degree at least $d/9$, by Chernoff bounds,
the probability that there is a vertex in $N'_{m+1}$ having
fewer than $d/36$ neighbours in either $P$ or $Q$ is at most
$2 \cdot (kd^2/729) \cdot e^{-d/72} < 1$, since $d \ge 300k \ge 600$.
Therefore we can choose a partition $N'_m = P \cup Q$ so that every vertex in $N'_{m+1}$
has at least $d/36$ neighbours in $P$  and at least $d/36$ neighbours in $Q$.

Consider the graph $G^*=G[P \cup N'_{m+1}]$ and its bipartite subgraph $B^*$ with parts
$P$ and $N'_{m+1}$, which has minimum degree $d^* \geq d/36 > 9k/2$. Fix any vertex
$x^*$ in $P$ and let $N^*_i$ denote the vertices at distance $i$ from $x^*$ in $B^*$.
By Lemma \ref{consecutive-paths}, we have some $m^* \geq 1$
such that $G^*$ contains cycles of length $\ell$ for all $3 \le \ell \le |N^*_{m^*}|/k$,
where $|N^*_{m^*}| \ge 2^{m^*-1}d^* \geq 2^{m^*}d/72$. We also have
a vertex $y$ such that there is a path between $x^*$ and $y$ in $G^*$
of length $\ell$, for all $\ell$ with $m^* \le \ell \le m^*+d/81-2 \leq m^*+4d^*/9-2$.
We let $y^*$ be either equal to $y$ if $y \in P$ or a neighbour of $y$
in $Q$ if $y \in N'_{m+1}$. In either case we have $y^* \in N'_m$
and there are paths between $x^*$ and $y^*$ in the bipartite subgraph $B_m$
of $G$ with parts $N_m$ and $N_{m+1}$ having any length $\ell$
with $m^*+1 \le \ell \le m^* + d/81 - 2$. Also, $x^*$ and $y^*$ both belong
to $N_m$, so are joined by a path $W$ of some length $\ell_W$ with $2 \le \ell_W \le 2m$,
where all internal vertices of $W$ lie in sets $N_i$ with $i<m$.
Combining $W$ with paths between $x^*$ and $y^*$ in $B_m$ gives cycles of
any length $\ell$ with $2m+m^*+1 \le \ell \le m^* + d/81$.
We already saw that $G$ contains cycles of length $\ell$ for all
$3 \le \ell \le \max\{|N_m|,|N^*_{m^*}|\}/k$.
Since $d \ge 300k$, we have $|N_m|/k \ge 2^{m-2}d/k>4m$ and $|N^*_{m^*}|/k \ge 2^{m^*}d/(72k)>4m^*$.
Therefore $\max\{|N_m|,|N^*_{m^*}|\}/k \ge 2m+m^*+1$, so $G$ contains cycles of
length $\ell$ for all $3 \le \ell \le d/81$. \qed

Next we need another lemma.

\begin{lemma}\label{paths+cycle}
Suppose $G$ is a graph with independence number $\alpha(G) \le k$
and $V(G)$ is partitioned into two parts $A$ and $B$ such that
\begin{itemize}
\item[(i)] $G[A]$ is Hamiltonian,
\item[(ii)] either $|B| \ge (9k+1)k+1$ or $G[B]$ has minimum degree at least $9k+1$,
and
\item[(iii)] every vertex in $B$ has at least $2$ neighbours in $A$.
\end{itemize}
Then $G$ contains a cycle of length $\ell$ for any
$2k+1+\lfloor\log_2(2k+1)\rfloor \le \ell \le |A|/2$.
\end{lemma}

\nib{Proof.}
First we note that $G[B]$ has an induced subgraph $H$
with minimum degree $d \ge 9k+1$ and at most $(9k+1)k+1$ vertices.
Indeed, if $|B| \le (9k+1)k+1$ just take $H=G[B]$,
otherwise apply Proposition \ref{indep-mindeg} to $G[B]$.
By Proposition \ref{bip}, $H$ contains a bipartite subgraph $H'$ with minimum degree $d/2> 9k/2$.
Applying Lemma \ref{consecutive-paths} to $H$ and $H'$, we find vertices $x$ and $y$ and a number
$m \le \log_2 \left( \frac{|V(H)|+d/2-1}{d/2} \right) \le \log_2(2k+1)$,
such that there is a path between $x$ and $y$ in $H$
of length $t$, for any $t$ with $m \le t \le m+2k-2 \leq m+\frac{4}{9}(d/2)-2$.
Since every $v \in S$ has at least $2$ neighbours in $A$ we can
choose neighbours $a$ of $x$ and $b$ of $y$ in $A$ with $a \ne b$.
Let $P$ be the path in $G[A]$ joining $a$ and $b$ obtained by taking the
longer arc of the Hamilton cycle, so that $P$ has length at least $|A|/2$.
We construct a cycle of any length $\ell$ with $2k+1+\log_2(2k+1) \le \ell \le |A|/2$ as follows.
Since $q=\ell-m-2k \geq 1$ we can apply Proposition \ref{shorten} to replace $P$ by a path $P'$
in $G[A]$ between $a$ and $b$ of some length $\ell'$ with $q  \leq \ell' \leq q+2k-2$.
Then $m \leq \ell-2-\ell' \leq m+2k-2$, so we can complete $P'$ to a cycle
of length $\ell$ by adding the edges $ax$, $by$ and a path in $H$ of length $\ell-2-\ell'$
between $x$ and $y$. \qed

Using this lemma we prove that if $G$ is a Hamiltonian graph on $n \ge 150k^3$ vertices with
independence number $\alpha(G) \le k$ then $G$ is pancyclic.

\noindent
{\bf Proof of Theorem \ref{pan-n}.}\,
Starting from the graph $G=G_n$ we construct a sequence of subgraphs $G_n, G_{n-1}, \cdots, G_{n-20k^2}$,
where $G_i$ is a Hamiltonian graph on $i$ vertices. Also, for each removed
vertex $v \in V(G) \sm V(G_i)$ we maintain a set of $2$ neighbours
$\{a_v,b_v\}$ of $v$ which we never delete, i.e.\ they will appear in each subgraph of the sequence.
To achieve this, consider the graph $G_i$, let $C_i$ be a Hamilton cycle in $G_i$
and let $N_i = \bigcup \big\{\{a_v,b_v\}: v \in V(G) \sm V(G_i) \big\}$.
We claim that we can choose $s=k^2+k+1$ disjoint $2k$-intervals in $C_i$ that avoid $N_i$.
To see this, consider the partition of $C_i$ into intervals defined by consecutive points in $N_i$.
We disregard $N_i$ and any intervals of length less than $2k$, then note that we can cover at least
half of the remaining points by disjoint $2k$-intervals. Since $|N_i| \le 40k^2$ and $n \ge 150k^3$
the number of $2k$-intervals thus obtained is at least $\frac{n-20k^2-40k^2(2k+1)}{2(2k+1)} > s$.
Now, as in the proof of Lemma \ref{erdos}, we can find a good jump $J$ in one of these intervals,
and use it to construct a cycle of length $i-1$. Furthermore, the vertex $v$ removed in this step
has at least $k+1$ neighbours in $V(C_i)$, since it belongs to the good jump $J$,
so we can choose any $2$ of these to be $a_v$ and $b_v$.

This sequence terminates with a Hamiltonian graph $G'=G_{n-20k^2}$ and a set
$S = V(G) \sm V(G')$ of size $20k^2>(9k+1)k+1$ such that every $v \in S$ has at least
$2$ neighbours in $V(G')$. By Lemma \ref{paths+cycle}
$G$ contains a cycle of length $\ell$ for any $2k+1+\log_2(2k+1) \le \ell \le |V(G')|/2=n/2-10k^2$.
To get cycles of length $\ell$ with $n/2-10k^2 \le \ell \le n$
we can just repeatedly apply Lemma \ref{erdos} starting from $G$.
To obtain the short cycles, note that $n \ge 150k^3 \ge (300k-1)k+1$, since $k \ge 2$,
so by Proposition \ref{indep-mindeg} $G$ has an induced subgraph $G^*$ with minimum degree
$d \ge 300k-1$. Since $(300k-1)/81 \ge 2k+1+\log_2(2k+1)$,
Theorem \ref{short-cycles} implies that $G^*$ contains cycles of length $\ell$
for all $3 \le \ell \le 2k+1+\log_2(2k+1)$. Therefore $G$ is pancyclic. \qed

Next we need the following lemma, which provides the long cycles
needed in the proof of Theorem \ref{pan-mindeg}.

\begin{lemma} \label{long}
Suppose $k \ge 3$ and $G$ is a Hamiltonian graph on $n \le 150k^3$ vertices
with minimum degree $\delta(G) \ge 600k$ and
independence number $\alpha(G) \le k$. Then $G$ contains
a cycle of length $\ell$ for any $n/12 \le \ell \le n$.
\end{lemma}

\nib{Proof.}
Consider a partition of the vertices of $G$ into sets $X$ and $Y$
where every vertex is placed randomly and independently in $X$ with probability $1/24$
or in $Y$ with probability $23/24$. By Chernoff bounds, the probability
that there is a vertex with less than $25k/2$ neighbours in $X$ is at most
$ne^{-25k/8}$ and the probability that $X$ has size more than $n/16$ is at most $e^{-n/384}$.
Since $k \ge 3$ and $600k \le n \le 150k^3$, both these probabilities are less than $0.4$,
so we can choose such a partition in which $|X| \le n/16$
and every vertex has at least $25k/2$ neighbours in $X$.
Starting from $G_n=G$ we construct a sequence of subgraphs
$G_n, G_{n-1}, \cdots, G_{n/12}$, where $G_i$ is a Hamiltonian graph
with $|V(G_i)|=i$ and $X \sub V(G_i)$.
To achieve this, suppose $n/12<i\le n$ and $C_i$ is a Hamilton cycle in $G_i$.
We claim that at least $3i/4$ of the vertices of $G_i$ are internal vertices
in some jump of length at most $8k$ in $C_i$. For if this is false, then by averaging
we could find an interval $I$ of length $8k$ and a set $S \sub I$ of size $2k+1$ such that
no vertex in $S$ is an internal vertex of a jump of length at most $8k$.
Consider a subset $S'$ of $S$ of size $k+1$ formed by taking every other vertex
(i.e.\ the first, the third, $\dots$, the $(2k+1)^{\text{st}}$).
Since $\alpha(G) \le k$ there must be an edge within $S'$.
This edge forms a jump of length at most $8k$ with at least one internal vertex in $S$,
giving a contradiction which proves the claim.
Since $i > n/12$ and $|X| \le n/16$ we have $|X| < 3i/4$. Thus we can choose a vertex $y \in Y$
and a jump $J$ of length at most $8k$ so that $y$ belongs to $J^o$ (the set of internal vertices of $J$).
Since every vertex has at least $25k/2$ neighbours in $X$ and $X \sub V(G_i)$,
every vertex in $J^o$ has at least $25k/2-8k>k+1$ neighbours in $V(G_i) \sm J^o$.
We replace the portion of $C_i$ traversing $J$ with the edge $\partial J$
then we use Lemma \ref{add} to put back the vertices of $J^o$ one by one,
until only $y$ has not been put back. Then $G_{i-1} = G_i \sm y$ is
Hamiltonian with $|V(G_{i-1})|=i-1$ and $X \sub V(G_{i-1})$, as required. \qed

Finally we give the proof of our third theorem, which states
that if $G$ is a Hamiltonian graph with minimum degree $\delta(G) \ge 600k$ and
independence number $\alpha(G) \le k$ then $G$ is pancyclic.

\noindent
{\bf Proof of Theorem \ref{pan-mindeg}.}\,
Let $n$ be the number of vertices in $G$. If $n \ge 150k^3$ then we are done
by Theorem \ref{pan-n} (without even using the minimum degree assumption)
so we can suppose that $n < 150k^3$. We also have $n > \delta(G) \ge 600k$, so $k \ge 3$.
Applying Theorem \ref{short-cycles} we see that $G$ has a cycle of length $\ell$
for all $3 \le \ell \le 7k <\delta(G)/81$. Also, by Lemma \ref{long} we have a cycle
of length $\ell$ for any $n/12 \le \ell \le n$.
For the remaining intermediate cycle lengths we consider a partition of the vertices
into two sets $X$ and $Y$, where vertices are randomly and independently placed
in $X$ with probability $1/2$ or in $Y$ with probability $1/2$.
By Chernoff bounds we can choose this partition so that $|X|, |Y| \ge n/3$
and each vertex has at least $200k$ neighbours in $X$ and at least $200k$ neighbours in $Y$.
Let $n'$ be the smallest number such that there is a subgraph $G'$ of $G$
on $n'$ vertices such that $G'$ is Hamiltonian and $X \sub V(G')$.
Let $C'$ be a Hamilton cycle in $G'$ and write $D = V(G) \sm V(G')$.
Since $X \sub V(G')$ we have $n' \ge n/3$ and $D \sub Y$.
First we dispose of the case when $|V(G') \cap Y| \le 4k$.
Since every vertex has degree at least $200k$ in $Y$,
the restriction of $G$ to $D=Y \sm V(G')$ has minimum degree at least $196k$.
Every vertex has at least $200k > 1$ neighbours in $X \sub V(G')$,
so applying Lemma \ref{paths+cycle} with $A = V(G')$ and $B=D$,
we obtain a cycle of length $\ell$ for any $2k+1+\log_2(2k+1) \le \ell \le n/6$.

Now we can suppose that $|V(G') \cap Y| \ge 4k+1$.
We can choose an interval $I$ of $C'$ that contains exactly $2k+1$ vertices of $Y$
and has length at most $n'/2$. Then we consider every other vertex of $Y$ in $I$
to obtain a set of size $k+1$, which must contain an edge, since $\alpha(G)\le k$.
This gives a jump $J$ of length at most $n'/2$ such that $1 \le |Y \cap J^o| \le 2k-1$.
Fix $y_0 \in Y \cap J^o$. We replace the portion of $C'$ traversing $J$ by $\partial J$,
and then use Lemma \ref{add} to put back vertices of $J^o \sm \{y_0\}$ one by one,
while we can find such a vertex with at least $k+1$ neighbours in the current cycle.
By minimality of $n'$ this process terminates before all vertices of $J^o \sm \{y_0\}$
have been replaced. Thus we obtain a non-empty subset $S$ of $J^o \sm \{y_0\}$ such that
$G''=G' \sm (S \cup \{y_0\})$ is Hamiltonian and every vertex in $S$ has at most $k$ neighbours in $V(G'')$.
Since $J$ has length at most $n'/2$ and $n' \geq n/3$ we have $n''=|V(G'')| \ge n/6$.
Also, $V(G'')$ contains $X \sm S$ and every vertex of $S$
has at least $200k$ neighbours in $X$, of which at most $k$ are in $G''$,
so the restriction $G[S]$ has minimum degree $d \ge 199k$.
Choose a vertex $x \in S$ that is adjacent to a vertex $a$ of $G''$. Such an $x$ exists
since $G'=G''\cup S \cup \{y_0\}$ is Hamiltonian, and in particular $2$-connected.
By Proposition \ref{bip} we can choose a bipartite subgraph $B$ of $G[S]$
with minimum degree at least $d/2$.
Applying Lemma \ref{consecutive-paths} to $G[S]$ and $B$
we obtain a vertex $y \in S$ and a number $m$
such that there is a path between $x$ and $y$ in $G[S]$ of length $\ell$,
for any $\ell$ with $m \le \ell \le m+4k$ (say, since $2d/9-2 > 4k$),
where $m \le \log_2 \left( \frac{n'/2+d/2-1}{d/2} \right) < \log_2 (150k^2/199+1) < 2\log_2 k$.

Let $C''$ be a Hamilton cycle in $G''$. If $y$ has a neighbour $b \ne a$ in $G''$
then we complete the argument as before. We take $P$ to be the longer arc of $C''$
between $a$ and $b$, so that $P$ has length at least $n''/2 \ge n/12$. Then we construct
a cycle of any length $\ell$ with $2k+1+2\log_2 k \le \ell \le n/12$ as follows. Since
$q = \ell-m-2k \ge 1$ we can apply Proposition \ref{shorten} to replace $P$ by a path $P'$ in $G''$
between $a$ and $b$ with some length $\ell'$ with $q \le \ell' \le q+2k-2$.
Then $m \le \ell-\ell'-2 \le m+2k$, so we can complete $P'$ to a cycle of length $\ell$
by adding the edges $ax$, $by$ and a path in $G[S]$ of length $\ell-\ell'-2$ between $x$ and $y$.
Now suppose that $y$ does not have a neighbour $b \ne a$ in $G''$.
We will repeatedly use the following fact.

\noindent $(\star)$ \ \
Any vertex $z$ with at most one neighbour in $G''$ has at least $40k$ neighbours in $D$.

The proof of $(\star)$ is immediate from that fact that $z$ has
at least $200k$ neighbours in $Y$, but at most $|Y \cap J|+1 \le 2k+2$
of these are in $G'$, so $z$ easily has at least $40k$ neighbours in $D$.
Applying $(\star)$ to $z=y$ we can choose a neighbour $y'$ of $y$ in $D$.
Let $Z$ be the connected component of $G[D]$ containing $y'$.
If $Z$ has an induced subgraph $Z'$ with minimum degree at least $20k$ then applying
Lemma \ref{paths+cycle} with $A = V(G')$ and $B=V(Z')$
gives a cycle of length $\ell$ for any $2k+1+\log_2(2k+1) \le \ell \le n/6$.

Now suppose that $Z$ does not have any induced subgraph $Z'$ with minimum degree at least $20k$.
We claim that there is a path of length at most $k$ in $Z$ from $y'$
to a vertex $z$ in $Z$ with at least $2$ neighbours in $G''$.
To see this note first that $|V(Z)| \le 20k^2$ by Proposition \ref{indep-mindeg}
and $Z$ contains at least one vertex with at least $2$ neighbours in $G''$ by $(\star)$.
Now consider a breadth first search tree $T$ in $Z$ starting from $y'$,
and for $i \ge 0$ let $N_i$ be the set of vertices at distance $i$ from $y'$
and let $Z_i$ be the restriction of $Z$ to $\cup_{j=0}^i N_j$.
If every vertex in $Z_i$ has at most one neighbour in $G''$ then by $(\star)$
it has at least $40k$ neighbours in $D$. On the other hand, we assumed that
the minimum degree in $G[Z_i]$ is less than $20k$. Therefore $Z_i$ contains a vertex $z$
with at least $40k-20k=20k$ neighbours in $D \sm
V(Z_i)$, implying $|N_{i+1}| \ge 20k$. Since $|V(Z)| \le 20k^2$ it follows that $Z_k$
contains a vertex $z$ with at least $2$ neighbours in $G''$, as
claimed.

Suppose that $z \in N_{i-1}$, where $1 \le i \le k+1$.
Combining the paths between $x$ and $y$ in $G[S]$ with the edge $yy'$
and the path of length $i-1$ in $T$ from $y'$ to $z$
we obtain paths between $x$ and $z$ in $G \sm V(G'')$ of any length $\ell$
with $m+i \le \ell \le m+4k+i$, where we recall that $m < 2\log_2 k$.
Choose a neighbour $b \ne a$ of $z$ in $G''$ and let $P$ be the longer arc of $C''$
between $a$ and $b$, so that $P$ has length at least $n''/2 \ge n/12$. (See Figure 1.)
Now we construct a cycle of any length $\ell$ with $3k+2+2\log_2 k \le \ell \le n/12$ as follows.
Since $q=\ell-m-2k-i \ge 1$ we can apply Proposition \ref{shorten} to replace $P$ by a path
$P'$ in $G''$ between $a$ and $b$ with some length $\ell'$ with $q \le \ell' \le q+2k-2$.
Then $m+i \le \ell-\ell'-2 \le m+2k+i$, so we can complete $P'$ to a cycle of length $\ell$
by adding the edges $ax$, $bz$ and a path in $G[S]$ of length $\ell-\ell'-2$ between $x$ and $z$.

\begin{figure}
\begin{center}
\includegraphics[scale=0.8]{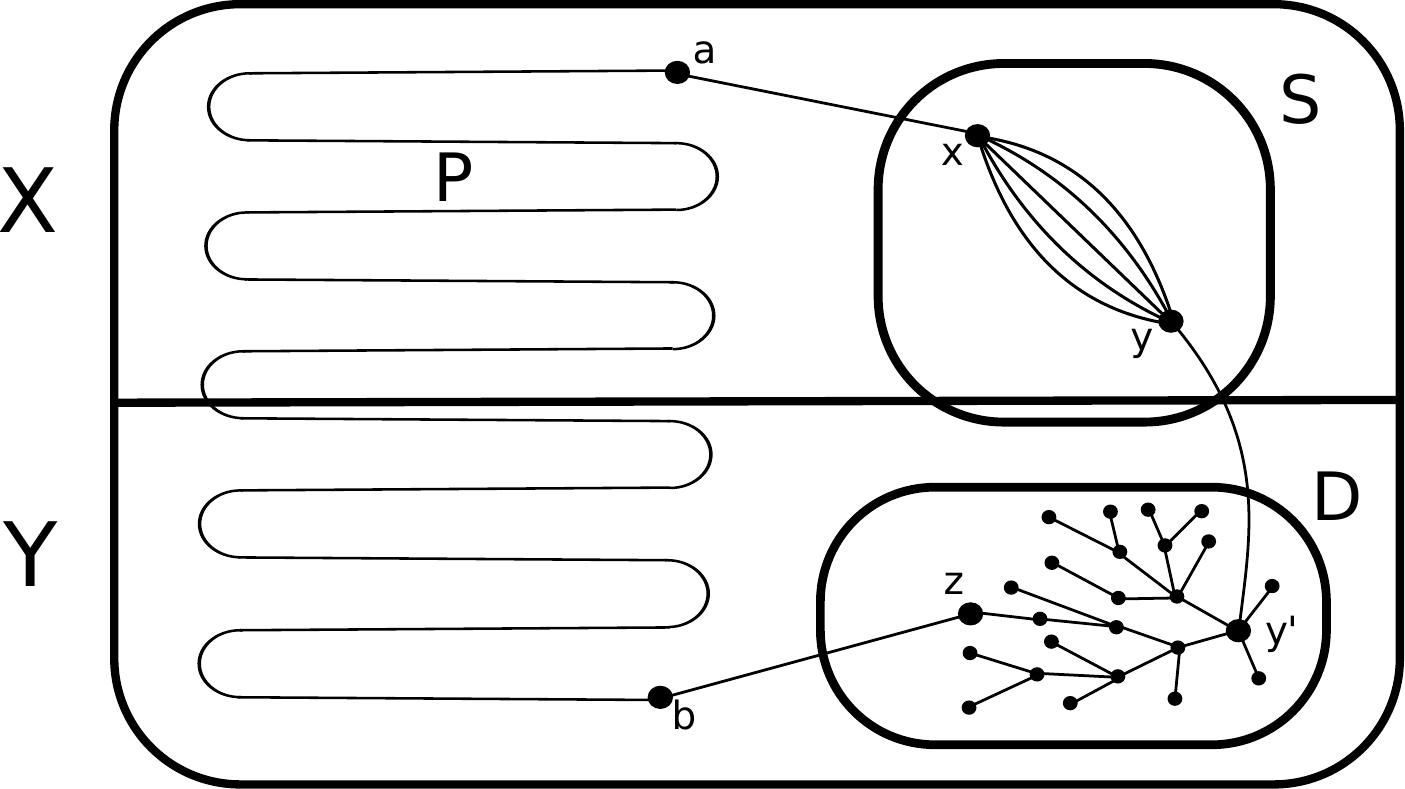}
\caption{Constructing cycles of intermediate length}
\end{center}
\end{figure}

Since $2k+1+2\log_2k<3k+2+2\log_2k<7k$, in all cases we find cycles of length $\ell$ for $7k \le \ell \le n/12$.
Recall that we also have cycles of length $\ell$ when $3 \le \ell \le 7k$
and when $n/12 \le \ell \le n$. This implies pancyclicity of $G$. \qed

\section{Concluding remarks}

We have answered the question of Jackson and Ordaz up to a constant factor.
Obviously it would be nice to obtain the exact bound, but perhaps one should
first attempt to prove an asymptotic version, i.e.\ that if
$\kappa(G) \ge (1+o(1))\alpha(G)$ then $G$ is pancyclic.
Also, it would be interesting to give the correct order of magnitude for
the minimum number $n$ of vertices such that any Hamiltonian graph $G$ on $n$
vertices with $\alpha(G)=k$ is pancyclic. We proved that this holds if $n=\Omega(k^3)$,
but it probably can be reduced to $n=\Omega(k^2)$. One way to attack this problem is
to improve the estimate in Lemma \ref{erdos}, which says that any Hamiltonian graph
with independence number $k$ and $n =\Omega(k^3)$ vertices contains a cycle of length $n-1$.
It would be extremely interesting to determine the correct dependence of $n$ on $k$
for this problem of just removing one vertex. Even the following question remains open.

\nib{Question.} Is there an absolute constant $C$ such that any Hamiltonian graph with
independence number $k$ and $n \geq Ck$ vertices contains a cycle of length $n-1$?

A positive answer would be tight up to a constant factor (clearly) and
in combination with Proposition \ref{indep-mindeg} and Theorem \ref{short-cycles}
would immediately imply that a quadratic dependence of $n$ on $k$ is sufficient
for Hamiltonicity to imply pancyclicity.


\begin{thebibliography}{99}

\bibitem{AS} N. Alon and J. Spencer,
{\em The probabilistic method}, second edition,
Wiley, New York, 2000.

\bibitem{AFG} D. Amar, I. Fournier and A. Germa,
Pancyclism in Chv\'atal-Erd\H{o}s's graphs,
{\em Graphs Combin.} {\bf 7} (1991), 101--112.

\bibitem{B1} J. A. Bondy, Pancyclic graphs I,
{\em J. Combin. Theory Ser. B} {\bf 11} (1971), 80--84.

\bibitem{B2} J. A. Bondy, Pancyclic graphs: recent results,
Infinite and finite sets, {\em Colloq. Math. Soc. J\'anos Bolyai},
181--187, Keszthely, Hungary, 1973.

\bibitem{B3} J. A. Bondy,
Basic graph theory: paths and circuits, Handbook of Combinatorics,
3--110, North-Holland, Amsterdam, 1995.

\bibitem{CE} V. Chv\'atal and P. Erd\H{o}s,
A note on Hamiltonian circuits,
{\em Discrete Math.} {\bf 2} (1972), 111--113.

\bibitem{D} G. A. Dirac, Some theorems on abstract graphs,
{\em Proc. London Math. Soc.} {\bf 2} (1952), 69--81.

\bibitem{E} P. Erd\H{o}s, Some problems in graph theory,
Hypergraph Seminar (Ohio State Univ., Columbus, Ohio, 1972),
{\em Lecture Notes in Math.} {\bf 411}, 187--190, Springer, Berlin, 1974.

\bibitem{EFRS} P. Erd\H{o}s, R. J. Faudree, C. C. Rousseau and R. H. Schelp,
On cycle--complete graph Ramsey numbers,
{\em J. Graph Theory} {\bf 2} (1978), 53--64.

\bibitem{FLMW}
E. Flandrin, H. Li, A. Marczyk and M. Wozniak,
A note on pancyclism of highly connected graphs,
{\em Disc. Math.} {\bf 286} (2004), 57--60.

\bibitem{G1} R. J. Gould,
Updating the hamiltonian problem - a survey,
{\em J. Graph Theory} {\bf 15} (1991), 121--157.

\bibitem{G2} R. J. Gould,
Advances on the hamiltonian problem - a survey,
{\em Graphs Combin.} {\bf 19} (2003), 7--52.

\bibitem{GHS} R. J. Gould, P. E. Haxell and A. D. Scott,
A note on cycle lengths in graphs,
{\em Graphs Combin.} {\bf 18} (2002), 491--498.

\bibitem{JO} B. Jackson and O. Ordaz,
Chv\'atal-Erd\H{o}s conditions for paths and cycles in graphs and digraphs: a survey,
{\em Disc. Math.} {\bf 84} (1990), 241--254.

\bibitem{L} D. Lou, The Chv\'atal-Erd\H{o}s condition for cycles in triangle-free graphs,
{\em Disc. Math.} {\bf 152} (1996), 253--257.

\bibitem{NS} V. Nikiforov and R. H. Schelp,
Cycle lengths in graphs with large minimum degree,
{\em J. Graph Theory} {\bf 52} (2006), 157--170.

\bibitem{SV} B. Sudakov and J. Verstra\"ete,
Cycle lengths in sparse graphs,
{\em Combinatorica} {\bf 28} (2008), 357--372.

\bibitem{V}  J. Verstra\"ete,
On arithmetic progressions of cycle lengths in graphs,
{\em Combin. Probab. Comput.} {\bf 9} (2000), 369--373.

\end{thebibliography}
\end{document}